\documentclass{article}
\usepackage[all]{xy}
\usepackage{graphicx}
\usepackage[margin=2.95cm]{geometry}
\usepackage{amssymb}
\usepackage{amsmath}
\usepackage{latexsym} 
\newtheorem{thm}{Theorem}[section]

\newtheorem{prop}[thm]{Proposition} 
\newtheorem{cor}[thm]{Corollary} \newtheorem{dfn}[thm]{Definition}

 \newtheorem{rmk}[thm]{Remark}
\newtheorem{ex}[thm]{Example} 

\newcommand{\pf}{\noindent{\bf Proof.}\ }

\newcommand{\integers}{{\mathbb Z}}

\newcommand{\Hom}{{\rm Hom}}

\newcommand{\im}{{\rm Im \,}}
\newcommand{\dom}{{\rm Dom \,}}
\newcommand{\ind}{{\rm Indet \,}}

\renewcommand{\ker}{{\rm Ker\, }}

\newcommand{\xbar}{{\overline{X}}}
\newcommand{\ybar}{{\overline{Y}}}
\newcommand{\zbar}{{\overline{Z}}}

\newcommand{\cald}{{\cal D}}
\newcommand{\cale}{{\cal E}}

\newcommand{\boldb}{{\mathbf B}}
\newcommand{\boldc}{{\mathbf C}}

\newcommand{\qed}{\begin{flushright} $\Box$\ \ \ \ \ \end{flushright}}

\newcommand{\from}{\leftarrow}

\newcommand{\rel}{\mathbf{REL}}

\newcommand{\lrel}{\mathbf{LREL}}

\newcommand{\clrel}{\mathbf{CLREL}}
\newcommand{\ilrel}{\mathbf{ILREL}}
\newcommand{\slrel}{\mathbf{SLREL}}

\newcommand{\unitob}{\mathbf{1}}

\newcommand{\arrows}{\,\lower1pt\hbox{$\longrightarrow$}\hskip-.24in\raise2pt
             \hbox{$\longrightarrow$}\,}
\newcommand{\twoheadleftarrowtail}{\twoheadleftarrow \mkern-15mu \leftarrowtail}

\newcommand{\xfy}{X \stackrel{f}{\longleftarrow}Y}

\newcommand{\ygz}{Y \stackrel{g}{\longleftarrow}Z}

\newcommand{\xfygz}{{X \stackrel{f}{\longleftarrow}Y \stackrel{g}{\longleftarrow}Z}}
\usepackage{url}

\title{Categories of (co)isotropic linear relations}
\author{Alan Weinstein\thanks{Research partially supported by UC Berkeley Committee on Research
\newline\mbox{~~~~}MSC2010 Subject Classification Number: 
 18B10 (Primary), 53D99 (Secondary)
\newline \mbox{~~~~}Keywords: linear relation, duality, symplectic  vector space,
coisotropic relation, Wehrheim-Woodward category}\\
Department of Mathematics\\University of California\\
Berkeley, CA 94720-3840 USA}
\date{}
\begin{document}
\maketitle

\begin{abstract}
In categories of linear relations between finite dimensional vector spaces, composition is well-behaved only
at pairs of relations satisfying transversality and monicity conditions.  A construction of Wehrheim and Woodward makes it possible to impose these conditions while retaining the structure of a category.   We analyze the resulting category in the case of all linear relations, as well as for (co)isotropic relations between symplectic vector spaces.  In each case, the Wehrheim-Woodward category is a central extension of the original category of relations by the endomorphisms of the unit object, which is a free submonoid with two generators in the additive monoid of pairs of nonnegative integers.
\end{abstract}

\section{Introduction}
Wehrheim and Woodward \cite{we-wo:functoriality} have introduced a category of symplectic manifolds allowing compositions between  canonical relations
(lagrangian submanifolds of products), such that the composition of relations $f$ and 
$g$ remains as  the usual set-theoretic composition when the pair $(f,g)$ satisfies local transversality and global embedding conditions.  For linear canonical relations between symplectic vector spaces, 
these two conditions are equivalent, and we have shown in \cite{li-we:selective} that the morphisms in the ``linear WW category" may be described as pairs consisting of a linear relation and a nonnegative integer which measures non-transversality and may be identified with an endomorphism of the zero-dimensional vector space.

In Poisson geometry, the most important relations (such as the graphs of Poisson maps) are the coisotropic ones, so it is interesting to study the WW category generated by such relations.  Here, the conditions of transversality and embedding are no longer equivalent, even in the linear case.  We will show that the description of a WW morphism in the linear coisotropic case requires {\it two} nonnegative integers.  For a pair $(A,B)$ of coisotropic subspaces representing an endomorphism of a point, these integers are the dimension of $A\cap B$ and the codimension of $A+B$.  (The numbers are equal in the lagrangian case.)

We will first study relations and then WW categories in the context of pure linear algebra, paying special attention to duality.
When we pass to symplectic linear algebra, as was the case for the classification of coisotropic pairs in \cite{lo-we:coisotropic},
it will be easier to deal directly with the isomorphic (via symplectic orthogonality) category of {\em isotropic} relations. 

We will end with some remarks on the question of extending our results from symplectic to presymplectic and Poisson vector spaces, as was done for the classification of subspace pairs in \cite{lo-we:coisotropic}.

In the original WW category, the endomorphisms of a point are equivalence classes of pairs of lagrangian \emph{submanifolds} in symplectic \emph{manifolds}; in that case, the equivalence relation remains quite mysterious, even for
endomorphisms of the one-point manifold represented by transversal lagrangian intersections in two points.  We hope that our analysis of the linear case will be a useful step toward dealing with the nonlinear case.

\section{Selective categories}
We begin by recalling some basic notions from \cite{li-we:selective}.  We refer to that paper for more details and further examples.

\begin{dfn}
\label{dfn-selective}
A {\bf selective category} is a category
 with a distinguished class of morphisms, called {\bf suave},
and a class of composable pairs of suave morphisms called {\bf
  congenial pairs}, such that:
\begin{enumerate}
\item Any identity morphism is suave.
\item If $f$ and $g$ are suave, $(f,g)$ is composable, and $f$ or $g$ is an identity morphism, then $(f,g)$ is congenial.
\item If $(f,g)$ is congenial, then $fg$ is suave.
\item If $f$ is a suave isomorphism, its inverse $f^{-1}$ is suave as well, 
and the pairs $(f,f^{-1})$ and $(f^{-1},f)$ are both congenial. 
\item If $(f,g,h)$ is a composable triple, then
  $(f,g)$ and $(fg,h)$ are congenial if and only if $(g,h)$
  and $(f,gh)$ are.   When these conditions hold, we call $(f,g,h)$ a
  {\bf congenial triple}. 
\end{enumerate}
\end{dfn}

It follows easily from the definition that any identity morphism is suave, and that notion of congeniality may be extended from pairs and triples to sequences of arbitrary length.

\begin{ex}
\label{ex-selectivesuave}
{\em 
We denote by  $\rel$ the category whose objects are sets and whose morphisms are
relations; i.e., the morphisms in $\rel(X,Y)$ are the subsets of
$X\times Y$.  
 
A useful selective structure in $\rel$ is the one in which all morphisms are suave, but only {\bf monic} pairs are congenial.   These are diagrams\footnote{We will usually denote morphisms by arrows pointing to the left, and we will therefore write 
$\Hom(X,Y)$ to denote the morphisms $X\from Y$  {\em to} $X$ {\em from} $Y$.}
 $\xfygz$ for which, whenever $(x,y)$ and $(x,y')$ belong to $f$,  
and $(y,z)$ and $(y',z)$ belong to $g$, then $y=y'$.  In other words, $(x,z)$ can belong to $f g$ via at most one element of the intermediate space $Y$.   The only condition in Definition  
\ref{dfn-selective} which takes a bit of work to check is the last one.   For this, it suffices to observe that congenial triples can be defined directly:
$(f,g,h)$ is congenial if and only if, whenever $(x,w)$ belongs to $fgh$, there is exactly one pair $(y,z)$ for which $(x,y) \in f$, $(y,z) \in g$, and $(z,w) \in h.$}
\end{ex}

\begin{dfn}
A {\bf selective functor} between selective categories is one which takes
congenial pairs to
congenial pairs.
\end{dfn}

Composition with identity morphisms shows that a selective functor takes
suave morphisms to suave morphisms. 
\bigskip

The Wehrheim-Woodward construction provides an embedding $\iota$ of the 
 the suave morphisms in
$\boldc$ to the morphisms in a category $WW(\boldc)$
in which the
composition of any {\em congenial} pair remains the same as the composition in $\boldc$.
Briefly stated (see \cite{li-we:selective} for details), the morphisms of $WW(\boldc)$ are equivalence classes of 
composable sequences in $C$, where two sequences are equivalent if one can be obtained from 
the other by moves in which a composable pair of adjacent entries in the sequence is replaced by its
composition, or vice versa.  In particular, identity morphisms may be freely inserted or removed.  We denote the equivalence class of $(f_1,\ldots,f_r)$ by  $[f_1,\ldots,f_r]$.

The map $\iota$ from the suave morphisms in $\boldc$
to $WW(\boldc)$ 
takes each suave morphism $f$ to $[f]$, the
equivalence class of the sequence having $f$ as its single entry.    $WW(\boldc)$ is then characterized by the following universal property.

\begin{prop}
\label{prop-universal}
Composition of the entries in composable sequences 
gives a well-defined
functor $\boldc \stackrel{c}{\from} WW(\boldc)$, 
$c([f_{1},\ldots,f_r]) =  f_{1}\cdots  f_r.$
 $WW(\boldc)$ has the universal property that 
any map from the suave morphisms of $\boldc$ to a category $\boldb$
which takes units to units and which takes congenial
compositions to compositions in $\boldb$ is of the form $b \iota$ 
for a unique functor $\boldb\stackrel{b}{\leftarrow}WW(\boldc).$ 
\end{prop}

We will refer to $c$ as the {\bf shadow} functor.   Composing it with the 
 inclusion $\iota$ of the suave morphisms in $\boldc$ gives the identity, so $\iota$ is injective, i.e., distinct suave morphisms cannot become equal when considered as WW-morphisms via $\iota$.  

\begin{rmk}
\emph{
Any selective functor between
selective categories induces a functor between their Wehrheim-Woodward
categories.  In particular, 
if a selective category $\boldc$ admits a transpose operation, i.e. an
involutive contravariant endofunctor $f\mapsto f^t$ which fixes objects and
takes congenial pairs to congenial pairs, then
this operation extends to $WW(\boldc)$, taking a sequence to 
 the same sequence with its order reversed and each entry replaced
by its transpose.   For instance, the operation which reverses the order of ordered pairs
is a transpose operation on $\rel$. }
\end{rmk}

\begin{rmk}
\label{rmk-subcategory}
{\em
If the suave morphisms in $\boldc$ form a subcategory $\boldc '$ and we
 declare every composable pair to be congenial, then 
$\boldc ' \stackrel{c}{\from} WW(\boldc)$ is an isomorphism of categories.
}
\end{rmk}

An important simplification result concerning selective categories requires the following finer structure.

\begin{dfn}
\label{dfn-highly}
A {\bf highly selective category} is a selective category provided with 
two subcategories of suave morphisms called {\bf reductions} and 
{\bf coreductions} such that:
\begin{enumerate}
\item Any suave isomorphism is a coreduction and a
  reduction.
\item If $(f,g)$ is a composable pair of suave morphisms, and if 
$f$ is a coreduction or $g$ is a reduction, then $(f,g)$ is
  congenial.
\item Any suave morphism $f$ may be factored as $gh$, where $g$ is a
  reduction, $h$ is a coreduction, and $(g,h)$ is congenial.
\end{enumerate}
\end{dfn}

It follows from the injectivity of $\iota$ that the subcategories of reductions and coreductions in $\boldc$ are mapped isomorphically to subcategories of $WW(\mathbf C)$, which we will again refer to as reductions and coreductions.   Since all identity morphisms are suave, these subcategories are wide; i.e., they contain all the objects.

We will indicate that a morphism is special by decorating its arrow:  $X \twoheadleftarrow Q$ is a reduction, and $Q \leftarrowtail Y$ is a coreduction.  A composition $X \leftarrow Q \leftarrow Y$ of suave morphisms is then congenial if at least one arrow is decorated at $Q$.

\begin{ex}
\label{ex-relhighlyselective}
{\em In the selective category of relations in Example \ref{ex-selectivesuave}, with the congenial pairs the monic ones, we obtain a highly selective structure by defining the reductions to be the single-valued (also called {\bf coinjective}) surjective relations and the coreductions the everywhere-defined (also called {\bf cosurjective}) injective ones.   }
\end{ex}

The following Theorem was proven in \cite{li-we:selective}, and in a special case already in \cite{we:note}.   It shows that arbitrary WW morphisms can be represented by sequences with just two entries.

\begin{thm}
\label{thm-two morphisms}
Let $\boldc$ be a highly selective category, and let 
$[f_1,\ldots,f_r]$ be a morphism in $WW(\boldc)$.  
Then there exist an object $Q$ and morphisms 
$A \in \Hom(X_0,Q)$ and $B \in \Hom(Q,X_r)$ in $\boldc$ such that  
$A$ is a reduction, $B$ is a coreduction, and 
$[f_1,\ldots,f_r]=[A,B].$
\end{thm}

A useful tool in the analysis of WW morphisms is the graph of a morphism.   This notion
applies when a category $\boldc$ has a rigid monoidal structure
compatible with its highly selective structure.   We recall that a monoidal structure on $\boldc$ consists of a ``tensor product'' bifunctor (often, but not always, denoted $\otimes$) from $\boldc \times \boldc$ to $\boldc$ and an object $\unitob$ which is a unit for this product.\footnote{For simplicity, we will assume that the associativity and unit axioms are satisfied strictly.}  The monoidal structure is rigid when each object $X$ has a dual $\xbar$ satisfying some identities which, in particular, give a bijective correspondence between $\Hom(X,Y)$ and $\Hom(X\otimes\ybar,\unitob)$.
The morphism $X\otimes \ybar \from \unitob$ corresponding to $\xfy$ in $\Hom(X,Y)$ is known as the {\bf graph} of $f$.   Compatibility with a selective structure means that the suave morphisms are closed under the operations and that a certain identity in the definition of rigidity involves only congenial compositions.  A key consequence of this compatibility is that a morphism $X\from Y$ is suave if and only if its graph $X\times \ybar \from \unitob$ is.

\section{Linear relations}
\label{sec-linear}

For simplicity, we will assume that all of our vector spaces are finite-dimensional.

\begin{dfn} The objects of the category $\lrel$ are the vector spaces over
a fixed ground field (which will usually remain unspecified), and
  the morphisms $\lrel(X,Y)$ are the linear subspaces of $X\times
  Y$.  Composition of morphisms is composition of relations.
\end{dfn}

The category $\lrel$ was studied in detail by Towber \cite{to:linear}
who gave an explicit classification of the endomorphisms up to conjugation by
automorphisms, as well as arbitrary morphisms up to ``left-right
equivalence''.   The endomorphism classification is also derivable
from the classification of subspace quadruples in \cite{ge-po:problems} since, as the 
authors  there observe, a relation $R\subseteq X\times X$ may be recovered from the quadruple of subspaces in $X \times X$ consisting of
$R$, $X\times 0_X$, $0_X\times X$, and the diagonal.

Since a given set may have more than one vector space structure, $\lrel$ is not a subcategory of
$\rel$, but the forgetful functor $\rel\from\lrel$ is faithful.  Via this functor, the special subcategories of
$\rel$ (injectives, surjectives, reductions, and coreductions)
and the monic pairs have their counterparts in $\lrel$.
Furthermore, we have the criteria that 
a linear relation $\xfy$ is injective if and only if its {\bf kernel},
the subspace $\ker f = f^t(0_X)$, 
is the zero subspace, and that $f$ is coinjective if and only if the kernel
$f(0_Y)$ of $f^t$ is zero.   We call this latter kernel the {\bf
  indeterminacy} of $f$ and denote it by $\ind f$.

Extending the well known fact that any linear map $\xfy$ may be
factored as $X\from \im f \from Y/\ker f \from Y$, we have the 
natural factorization of any linear relation $\xfy$ as the
strongly transversal composition 
\begin{equation}
\label{eq-factorlinearinthree}
X\leftarrowtail \im f/\ind f \twoheadleftarrowtail \dom
f/\ker f \twoheadleftarrow Y.
\end{equation}
This factorization leads to a proof of the fact stated at the end of
Section 1 of \cite{to:linear} that a linear relation is classified up
to isomorphism by the dimensions of the six spaces occurring in
(\ref{eq-factorlinearinthree}). 

 \subsection{Monicity}

A composable pair $(f,g)$ of linear relations
is monic 
if and only 
 $\ker f \cap \ind g = 0$, in particular if $f$ is injective or $g$
is coinjective.   Another criteron for 
monicity of the pair $(f,g)$
when $f\in \lrel(X,Y)$ and $g\in \lrel(Y,Z)$ is
that the subspaces $f\times g$ and $0_{X}\times \Delta_Y\times {0_Z}$ of
$X\times Y\times Y\times Z$ intersect in the zero subspace.  In general, the dimension of this intersection
will be called the {\bf defect} of the pair.

We may extend the definition and quantitative measure of monicity to composable sequences of arbitrary length.

\begin{dfn}
The dimension of
$$(f_1 \times f_2 \times \cdots \times f_r)\cap
(0_{X_0} \times
\Delta_{X_1}\times\cdots\times \Delta_{X_{r-1}}\times 0_{X_r})$$ will be called the {\bf excess} of the sequence
and denoted by $\cale(f_1,f_2,\ldots, f_r)$.  
\end{dfn}

The excess  $\cale(f_1,f_2,\ldots, f_r)$ may be interpreted as the dimension, for any $(x_0,x_r) $ in the composed relation $f_1\cdots f_r \subseteq( X_0\times X_r)$,  of the affine space of $r$-tuples
$$((x_0,x_1),(x_1,x_2),\ldots, (x_{r-2},x_{r-1}),(x_{r-1},x_r)) $$
with each $(x_{i-1},x_i) \in f_i)$.   We will refer to such $r$-tuples as {\bf trajectories} to $x_0$ from $x_r$.
We will use the term ``monic" by extension from the case $r=2$ for any composable sequence with zero excess, i.e. with at most one trajectory between any two endpoints.   Notice that any 1-entry sequence has zero excess.   

\begin{prop}
\label{prop-composedexcess}
Let $(f_1,f_2,\ldots, f_r)$ be a composable $r$-tuple of linear relations, and let $1\leq j \leq r-1$.
Then $\cale(f_1,f_2,\ldots, f_r) = \cale(f_1,f_2,\ldots, f_j) +\cale(f_1\cdots f_j,f_{j+1}\cdots f_r) + \cale(f_{j+1},f_{j+2},\ldots, f_r)$.
\end{prop}

\pf
We consider $\cale(f_1,f_2,\ldots, f_r)$  as the dimension of the (linear) space of trajectories to $0_{X_0}$ from $0_{X_r}$.   Each such trajectory can be decomposed into a pair of trajectories, $\tau'$ to $0_{X_0}$ from some $x_j \in X_j$, and $\tau'' $ to this $x_j$ from $0_{X_r}$.   The set $S_j$ of $x_j$ which occur in this way may be identified with the space of trajectories to $0_{X_0}$ from $0_{X_r}$ for the 
2-member sequence $(f_1\cdots f_j,f_{j+1}\cdots f_r)$; it therefore has dimension $\cale(f_1\cdots f_j,f_{j+1}\cdots f_r)$.   For each $x_j \in S_j$, the affine spaces $\tau'$ and $\tau''$ of trajectories which connect it to $0_{X_0}$ and $0_{X_1}$ have dimensions $\cale(f_1,f_2,\ldots, f_j)$ and $\cale(f_{j+1},f_{j+2},\ldots, f_r)$ respectively.  Adding the dimensions of these three affine spaces gives the result.
\qed

This leads immediately to:

\begin{cor}
\label{cor-highertransversal}
Let $(f_1,f_2,\ldots, f_r)$ be a composable $r$-tuple of linear relations, and
consider some implementation  of the composition $f_1 f_2 \cdots f_r$ 
by a sequence of pairwise compositions, i.e. a ``parenthesization.'' 
Then the excess of the composable $r$-tuple may be obtained
by summing the excesses of the pairwise compositions and  is therefore independent of the parenthesization.  In particular,  the $r$-tuple is monic
if and only if all of the pairwise compositions in the sequence are monic.  
\end{cor}

\subsection{Duality}
\label{subsec-duality}

The contravariant transposition functor on $\lrel$ which exchanges the entries in each ordered pair  is
compatible with the one on $\rel$ via the forgetful functor between
these two categories, but there is another contravariant functor on $\lrel$
which does not correspond to anything in $\rel$.  It extends the usual
duality for linear maps.

\begin{dfn}
For $f\in
\lrel(X,Y)$, the {\bf dual} $f^*\in \lrel(Y^*,X^*)$ is the relation 
$\{(\eta,\xi)\in Y^*\times X^* | \eta(y)=\xi(x) \mbox{~for all~}
(x,y)\in f \}.$
\end{dfn}

\begin{rmk}
\label{rmk-dualannihilator}
\emph{The dual $f^*$ may also be seen as the annihilator of $f$ with respect to
the bilinear pairing $$\langle(x,y),(\eta,\xi)\rangle =
\xi(x)-\eta(y)$$ between $X\times
Y$ and $Y^*\times X^*$.  When $X$
and $Y$ are reflexive, this pairing is nondegenerate, from which it
follows easily that  $f^{**}=f$ when we identify $X^{**}$ with $X$ and $Y^{**}$ with
$Y$.}  
\end{rmk}

That duality is a (contravariant) functor is not completely trivial to prove.  In fact, it
depends on the fact that any linear functional on a subspace of a
vector space can be extended to the whole space;
functoriality fails for modules over arbitrary rings exactly when
there are nonextendible functionals.

Proofs of the following result may be found in
 \cite{be-tu:relazioni} and  \cite{to:linear}; the latter applies in the
infinite-dimensional case as well.

\begin{prop}
\label{prop-functorialdual}
Duality is functorial; i.e., 
for  linear relations  $\xfygz$ , $(f g)^* = g^* f^*$.
\end{prop}

\begin{ex}
{\em We may also 
consider the category of linear relations in modules over any
  ring $R$.  Let $X$ be a submodule of $Y$ such that there is an
  element $\xi_0$ of
  $X^*$ which does not extend to $Y$.  
(For instance, with $R=\integers$, let $X=2 \integers$,
$Y=\integers$, and $\xi_0(n)=n/2$.)  
Let $Z$ be the zero module.
Then, if $f\in \lrel(X,Y)$ is
  the transpose of the inclusion of $X$ in $Y$, and $g\in \lrel(X,Y)$
  is the graph of the zero map, $(f g)^* \neq g^* f^*.$

In fact, $f g$ is the graph of the zero map $X\from Z$, and hence
$(f g)^*$ is the graph of the zero map $Z^* \from X^*$.  On the
other hand,  $g^*$ is the graph of the zero map $Z^*\from Y^*$,
while
$f^*$ is the transpose of the restriction map $X^*\from
Y^*.$  The composition $g^* f^*$ thus consists of those pairs
$(0_{Z^*},\xi)\in Z^* \times X^*$ for which there is an $\eta\in Y^*$ whose
restriction to $X$ is $\xi$.  Then $(0_{Z^*},\xi_0)$
belongs to $(f g)^*$ but not to $g^* f^*.$
 }
\end{ex}

The
duality functor has the following properties (some true for modules as
well as vector spaces).

\begin{prop}
For any linear relation $\xfy$ between vector spaces:
\begin{enumerate}
\item 
$(f^t)^* = (f^*)^t$;
\item
$\ker f^* =  (\im f)^\circ$;
\item
$\im f^* = (\ker f)^\circ$;
\item
$\ind f^* = (\dom f)^\circ$;
\item
$\dom f^* = (\ind f)^\circ$,
\end{enumerate}
where $^\circ$ denotes the annihilator of a subspace in the dual of the ambient space.

Consequently, the duality functor exchanges injectives with
cosurjectives, coinjectives with surjectives, and reductions with coreductions.
\end{prop}

\pf
The first two properties are immediate from the definitions.  For 3,
one uses the equation $\xi(x) = \eta(y)$ for $(x,y) \in f$ 
to define a functional $\xi$
on the image of $f$ if $f$ annihilates the kernel of $f$
and then extends $\xi$ to all of $X$.  
Transposing statements 2 and 3  yields statements 4 and 5.
\qed

We also have:

\begin{prop}
\label{prop-dualmonic}
Let $(f,g)$ be a composable pair of linear relations.  Then:
\begin{itemize}
\item
$(\ker f \cap \ind g)^\circ  = \dom g^* +\im f^* $, and
\item
$(\dom f +\im g) ^\circ = \ker g^*\cap \ind f^*  $
\end{itemize}
\end{prop}

\begin{rmk}
\label{rmk-pairing}
\emph{
All of the results in this section remain valid if $V^*$ is replaced by any space equipped with a nondegenerate pairing with $V$, such as $V$ itself when it carries a symplectic structure.}
\end{rmk}

\subsection{Transversality}
To define selective structures on categories of linear relations, we can of course use the condition of monicity, carried over from $\rel$, but compatibility with duality and other considerations (such as continuity) suggest that we impose
another condition as well, which does not have a counterpart in $\rel$.

\begin{dfn}
\label{dfn-transversal}
A pair $(f,g)$ of linear relations $\xfy$ and $\ygz$ is {\bf transversal} if 
the sum $f\times g +X \times
\Delta_Y \times Z$ is the entire space
$X\times Y\times Y \times Z$.  A pair which is both transversal and
monic is {\bf strongly transversal}.
\end{dfn}

The pair $(f,g)$ of linear relations  is 
transversal if and only if $\dom f + \im g$
is the entire intermediate space $Y$.  
In particular, $(f,g)$ is transversal
whenever $f$ is cosurjective or $g$ is surjective.

Recall that the composition of any sequence $(f_1,f_2,\ldots, f_r)$ of linear relations, with $f_j \in
\lrel(X_{j-1},X_j)$ may be obtained by taking the intersection of 
 $f_1 \times f_2 \times \cdots \times f_r$ with $X_0 \times
\Delta_{X_1}\times\cdots\times \Delta_{X_{r-1}}\times X_r$ in
 $X_0 \times X_1 \times\cdots
X_r$ and then projecting this along its intersection
with the kernel
$0_{X_0}\times
\Delta_{X_1}\times\cdots\times \Delta_{X_{r-1}}\times 0_{X_r}$ of the projection to
 $X_0 \times X_r$.

The codimension of $$ (f_1 \times f_2 \times \cdots \times f_r) + (X_0 \times
\Delta_{X_1}\times\cdots\times \Delta_{X_{r-1}}\times X_r)$$ in
 $X_0 \times X_1 \times\cdots
X_r$ will be called the {\bf defect} of the sequence and denoted by $\cald(f_1,f_2,\ldots, f_r)$. 
We will call the entire sequence transversal if its defect is zero; this extends the definition for 2-entry sequences.

Using properties of duality, we can now transfer statements about excess and monicity to defect and transversality.
 
\begin{prop}
\label{prop-defectexcess}
Let $(f_1,f_2,\ldots, f_r)$ be a composable $r$-tuple of linear relations.  \\Then
$\cale(f_1,f_2,\ldots, f_r)= \cald(f_r^*,\ldots,f_1^*)$, and 
$\cald(f_1,f_2,\ldots, f_r)=\cale(f_r^*,\ldots,f_1^*)$.

In particular, a composable pair $(f,g)$ is monic if and only if
$(g^*,f^*)$ is transversal, and vice versa.  Thus, the duality
functor preserves strong transversality.

\end{prop}
\pf
The excess $\cale(f_1,f_2,\ldots, f_r)$ is the dimension of the intersection
$$(f_1 \times f_2 \times \cdots \times f_r)\cap
(0_{X_0} \times
\Delta_{X_1}\times\cdots\times \Delta_{X_{r-1}}\times 0_{X_r}).$$ 

We now pair the spaces $$X_0\times X_1\times X_1\times X_2\cdots \times X_{r-1}\times X_{r-1}\times X_r$$
and $$X_r^*\times X_{r-1}^*\times X_{r-1}^*\times  X_{r-2}^* \times  \cdots \times X_{2}^*\times X_{1}^*
\times X_{1}^*\times X_0^*$$
by the rule
\begin{eqnarray*}\langle  (x_0, x_1',x_1, x_2',\ldots , x_{r-1}',x_{r-1},x_r'),
(\xi_r, \xi_{r-1}', \xi_{r-1}, \xi_{r-2}',  \ldots , \xi_{1}',
\xi_{1},\xi_0' )                 \rangle \\ 
=\xi_r (x_r') - \xi_{r-1}' (x_{r-1}) + \xi_{r-1} (x_{r-1}') - \xi_{r-2}' (x_{r-2}) +\cdots
       + \xi_1 (x_1') - \xi_{0}' (x_{0})
 \end{eqnarray*}

The annihilators of the spaces being intersected are now
$$f_r^* \times f_{r-1}^* \times \cdots \times f_1^* \mbox{~~and~~}
X_r^* \times
\Delta_{X_{r-1}}\times\cdots\times \Delta_{X_{1}}\times X_0^*$$
respectively, the codimension of whose sum is the defect $\cald(f_r^*,\ldots,f_1^*)$.
Since, for a nondegenerate pairing, the dimension of an intersection is the codimension of the sum of
the annihilators, our equation is proven.

The second equation in the statement of the proposition follows immediately from the first by reflexivity.
\qed

\begin{cor} 
\label{cor-composeddefect}
Let $(f_1,f_2,\ldots, f_r)$ be a composable $r$-tuple of linear relations,
and let $1\leq j \leq r-1$.
Then $\cald(f_1,f_2,\ldots, f_r) = \cald(f_1,f_2,\ldots, f_j) +\cald(f_1\cdots f_j,f_{j+1}\cdots f_r) + \cald(f_{j+1},f_{j+2},\ldots, f_r)$.
\end{cor}
\pf
Combine Propositions \ref{prop-composedexcess} and \ref{prop-defectexcess}. 
\qed

\section{Isotropic and coisotropic relations}
We now define two enlargements of the category $\slrel$ of symplectic vector spaces and canonical relations.
In the categories
$\ilrel$ and $\clrel$, the objects are again the symplectic vector spaces, and the morphisms $X\from Y$ are the subspaces of $X\times \ybar$ which are isotropic in the case of $\ilrel$ and coisotropic in the case of $\clrel$.   It is not hard to see that these collections of subspaces  are closed under composition and contain the diagonals.   The category $\slrel$ then sits naturally as a subcategory of $\ilrel$ and $\clrel$ and is the intersection of those two categories.

If $I$ is isotropic in $X$, then its annihilator with respect to the symplectic pairing is coisotropic, and vice versa.  Since the duality functor may be defined in terms of annihilators, we have a bijective functor between $\ilrel$ and $\clrel$ which is the identity on objects.  It is convenient to define this functor in contravariant form in such a way that it becomes the identity on objects, as well as on their intersection $\slrel$.  

Since a lagrangian relation is isomorphic to its dual, it follows from Proposition \ref{prop-dualmonic} that the defect and excess are equal in $\slrel$.  More generally, we have:

\begin{prop}
\label{prop-defectexcessinequal}
The excess of a composable sequence $(f_1,f_2,...,f_r)$ in $\ilrel$ is less than or equal to its defect.  In $\clrel$, the defect is less than or equal to the excess.  These are the only constraints on the defect and excess.
\end{prop}
\pf
The case of $\clrel$ follows from that of $\ilrel$ by Proposition \ref{prop-dualmonic}.  The defect $\cald(f_1,f_2,...,f_r)$ is
 the codimension of $$ (f_1 \times f_2 \times \cdots \times f_r) + (X_0 \times
\Delta_{X_1}\times\cdots\times \Delta_{X_{r-1}}\times X_r)$$ in
 $X_0 \times X_1 \times\cdots
X_r$.   The codimension of this space is equal to the dimension of its symplectic orthogonal
$$ (f_1^\omega \times f_2^\omega \times \cdots \times f_r^\omega) \cap (0_{X_0} \times
\Delta_{X_1}\times\cdots\times \Delta_{X_{r-1}}\times 0_{X_r}).$$   When the $f_j$ are isotropic, this space
contains
$$ (f_1\times f_2\times \cdots \times f_r) \cap (0_{X_0} \times
\Delta_{X_1}\times\cdots\times \Delta_{X_{r-1}}\times 0_{X_r}),$$ whose dimension is the 
excess of $(f_1,f_2,...,f_r)$.  Hence $\cale(f_1, f_2, \cdots, f_r) \leq \cald (f_1, f_2, \cdots, f_r).$

Simple examples show that there are no further constraints on $\cald$ and $\cale$.
\qed

\section{Selective structures and WW morphisms}

We showed in \cite{li-we:selective} that the morphisms in $WW(\slrel)$ may be identified with {\bf indexed canonical relations}, which are pairs $(L,k)$ such that $L$ is a canonical relation and $k$ is a nonnegative integer.  The composition 
$(L,k)(L',k')$ is defined, whenever $(L,L')$ is a composable pair, as  $(LL',k+k'+\cale(L,L'))$.
This makes $WW(\slrel)$ a central extension\footnote{Central extensions of categories are defined
in \cite{ne:lectures}.} of $\slrel$ by the endomorphisms of the unit object,  which may identified with the additive monoid $\integers_+$ of nonnegative integers and which are composed with other morphisms by the tensor product.  The excess $\cale$, equal here to the defect $\cald$, is the cocycle determing  the extension.

When we replace $\slrel$ by one of the categories $\lrel$, $\ilrel$, or $\clrel$, the excess and defect are no longer equal.
We will therefore use the highly selective structures on $\lrel$, $\ilrel$ and $\slrel$ in which all morphisms are suave, while the congenial pairs are only those which are  strongly transversal, i.e. both monic and transversal. The reductions are again the surjective and single-valued  relations, and the coreductions are the everywhere-defined and injective relations.

\begin{rmk}
\emph{The cotangent functor $T^*$ 
from $\lrel$ to $\slrel$ takes a vector space $V$ to the symplectic 
space $V\oplus V^*$ with the usual structure and a linear relation $V\from W$ to the canonical relation which is its conormal bundle with the components in $W^*$ multiplied by $-1$.\footnote{This is a very special case of the symplectic groupoid construction taking a coisotropic relation between Poisson manifolds to a canonical groupoid relation between their symplectic groupoids.}
For a composable pair $(f,g)$ of linear relations, the defect (equal to the excess) of the pair $(T^*f,T^*g)$ is the sum of the defect and excess of $(f,g)$.   
}
\end{rmk}

Since vanishing of the defect and excess are the defining conditions for congeniality, it is perhaps not surprising that it suffices to keep track of just these two numbers to parametrize the WW morphisms.  

\begin{thm}
\label{thm-doublyindexed}
If $\boldc$ is the category $\lrel$, $\ilrel$, or $\clrel$, any morphism in  $WW(\boldc)$ is 
determined by its shadow, its defect, and its excess.  These are independent, subject only to the inequalities in Proposition \ref{prop-defectexcessinequal}.   The category $WW(\boldc)$ is 
a central extension of $\boldc$ by the endomorphisms of the unit object.  In the case if $\lrel$, the map $(\cald,\cale)$ identifies these endomorphisms with the additive monoid $\integers_+ \times \integers_+$.   For $\ilrel$ or $\clrel$, we get  the submonoid defined by the inequalities 
$\cale \leq \cald$ or $\cald \leq \cale$.   The cocycle determining the central extension is $(\cald,\cale)$.
\end{thm}

\pf
We begin with the case of $\lrel$.
Any morphism in $WW(\lrel)$ is represented by a diagram
$$X \stackrel {f}{\twoheadleftarrow}Y\stackrel{g}{\leftarrowtail}Z$$ of linear relations.  The graph of this composition may be constructed from the graphs $\gamma_f$ and $\gamma_g$ and is represented by the diagram of linear relations
$$X \times Z \stackrel{X\times \Delta_Y\times Z}{\twoheadleftarrow}
X \times Y\times Y\times Z
                   \stackrel{(\gamma_f \times \gamma_g)\delta_\unitob}{\leftarrowtail}\unitob,$$
where $X\times \Delta_Y\times Z$ denotes the reduction whose domain is this subspace
and whose kernel is $0_X\times\Delta_Y\times 0_Z$, and $\delta_\unitob$ is the diagonal isomorphism $\unitob \times \unitob \from \unitob$.  

Analyzing this composition comes down to analyzing the triple $(X\times \Delta_Y\times Z, 0_X\times\Delta_Y\times 0_Z,f \times g)$ of subspaces in $X\times Y \times Y\times Z$, using the natural isomorphism between $X\times Z$ and
 $(X\times\Delta_Y\times Z)/ (0_X\times\Delta_Y\times 0_Z)$.  
Such triples have been classified in detail, for instance in \cite{br:endomorphisms} and
\cite{et:introduction}.  In this classification, any quadruple $(V,A,B,C)$ consisting of a vector space $V$ and its subspaces $A$, $B$, and $C$ is isomorphic to a direct sum of copies of ``elementary"  components $(V_i,A_i,B_i,C_i)$ of the following nine types.  $V$ is 1-dimensional in the first eight types and is 2-dimensional in the ninth.

\begin{description}
\item[$\tau_1$:] $A = 0$, $B=0$,  $C=0$
\item[$\tau_2$:] $A = 0$, $B=0$,  $C=V$
\item[$\tau_3$:] $A = 0$, $B=V$,  $C=0$
\item[$\tau_4$:] $A = V$, $B=0$,  $C=0$
\item[$\tau_5$:] $A = 0$, $B=V$,  $C=V$
\item[$\tau_6$:] $A = V$, $B=0$,  $C=V$
\item[$\tau_7$:] $A = V$, $B=V$,  $C=0$
\item[$\tau_8$:] $A = V$, $B=V$,  $C=V$
\item[$\tau_9$:] $A$, $B$, and $C$ are pairwise independent subspaces, any two of which sum to $V$.
\end{description}

Applying this decomposition to the case $V = X\times Y \times Y\times Z$, 
$A=X\times \Delta_Y \times Z$, 
$B=0_X\times \Delta_Y\times 0_Z,$ $C=f \times g$ and using the compatible decompositions of $A$ and $B$ to obtain a decomposition of $X\times Z=A/B$
 we obtain a decomposition of our original WW morphism as a direct sum of nine ``elementary morphisms", each of which may contribute to the shadow, defect, and excess of the morphism.    Since
 $B$ is contained in $A$, the types $\tau_3$, $\tau_5$, and $\tau_9$ cannot appear in the decomposition.
 
 The contributions of the remaining six types are as follows, with the dimension of $V_i$ denoted by $n_i$:

 \begin{center}
\begin{tabular}{ c|c|c|c|c } 
type & $\dim(A/B)$ & $\dim$(shadow) &defect & excess\\ 
 \hline
 $\tau_1$&0  &               0 &$n_1$&0\\ 
$\tau_2 $& 0 &               0 & 0&0\\
$\tau_4 $& $n_4$ &     0  &0&0\\ 
$\tau_6 $&  $n_6$& $n_6 $&0&0\\  
$\tau_7$ & $0$ &        0  &0&0\\  
$\tau_8$ & $0$ &       0  &0&$n_8$\\ 

 \hline
\end{tabular}
\end{center}

The components of type $\tau_1$ and $\tau_8$ are the graphs of endomorphisms of the unit object $\unitob$ which are determined by their defect and excess.  The components of type $\tau_2$ and $\tau_7$ may be ignored, since they give the trivial endomorphism of $\unitob$ as a morphism in $WW(\lrel)$. Finally, the components of types 
$\tau_4$ and $\tau_6$ are determined by the shadow of $(f,g)$, with dimension $n_6$, and the space $X\times Z$, with dimension $n_4 + n_6. $  It follows that the graph of our WW morphism, and thus the morphism itself, is determined by the shadow, the defect, and the excess.    Examples showing that all possibilities for these three items occur are easy to construct, and the composition law based on the cocycle $(\cald,\cale)$ follows from Proposition
 \ref{prop-composedexcess} and Corollary \ref{cor-composeddefect}.
\bigskip

We now go on to the case of $WW(\ilrel)$; the case of $WW(\clrel)$ is reduced to this one by taking symplectic orthogonals.  
We will modify the argument above, taking the extra symplectic data into account.  Since the kernel of a reduction is determined by the domain, being its symplectic orthogonal,  we may work with pairs of subspaces instead of triples.

Each morphism in $WW(\ilrel)$ is represented by a diagram
$$X \stackrel {f}{\twoheadleftarrow}Y\stackrel{g}{\leftarrowtail}Z$$ of isotropic relations.  As above,
the graph of this composition may be constructed from the graphs $\gamma_f$ and $\gamma_g$ and is represented by the diagram of isotropic relations
$$X \times \zbar \stackrel{X\times \Delta_Y\times \zbar}{\twoheadleftarrow}
X \times \ybar\times Y\times \zbar
                   \stackrel{(\gamma_f \times \gamma_g)\delta_\unitob}{\leftarrowtail}\unitob,$$
where $X\times \Delta_Y\times \zbar$ denotes the (lagrangian, hence isotropic) reduction by this coisotropic subspace, and $\delta_\unitob$ is the diagonal isomorphism $\unitob \times \unitob \from \unitob$.  

Analyzing this composition comes down to decomposing the pair $(X\times \Delta_Y\times \zbar,f \times g)$ of subspaces in $X\times\ybar\times Y\times \zbar$, using the natural isomorphism between $X\times \zbar$ and
 $(X\times\Delta_Y\times \zbar)/ (0_X\times\Delta_Y\times 0_\zbar)$.   In fact, it is more convenient to pass to the equivalent problem of decomposing the isotropic pair $(0_X\times \Delta_Y\times 0_\zbar,f\times g)$, since pairs of isotropic subspaces have already been classified explicitly in \cite{lo-we:coisotropic}.   In this classification, any triple $(V,A,B)$ consisting of a symplectic space $V$ and its isotropic subspaces $A$ and $B$ is isomorphic to a monoidal product of ``elementary"  components $(V_i,A_i,B_i)$ of the following five types:

\begin{description}
\item[$\tau_1$:] $A$ and $B$ are lagrangian subspaces, and $A = B$
\item[$\tau_2$:] $B = 0$ and $A$ is a lagrangian subspace
\item[$\tau_3$:] $A = 0$ and $B$ is a lagrangian subspace
\item[$\tau_4$:] $A$ and $B$ are lagrangian subspaces, and $A \cap B = 0$
\item[$\tau_5$:] $A = B = 0$
\end{description}

Applying this decomposition to the case $V = X\times\ybar\times Y\times\zbar$, $A=0_X\times \Delta_Y\times 0_\zbar,$ $B=f \times g$ and using the compatible decompositions of $A$ and $A^\omega$ to obtain a decomposition of $X\times\zbar=A^\omega/A$
 we obtain a decomposition of our original WW morphism as a monoidal product of five ``elementary morphisms", each of which may contribute to the shadow, defect, and excess of the morphism.    The contributions are now
as follows, with {\em half} the dimension of $V_i$ denoted by $n_i$.

\begin{center}
\begin{tabular}{ c|c|c|c|c } 
type & $\dim(A^\omega/A)$ & $\dim$(shadow) &defect & excess\\ 
 \hline
 $\tau_1$&0  &0 &$n_1$&$n_1$\\ 
$\tau_2 $& 0 &  0 &$n_2$&0\\
$\tau_3 $& $2n_3$ &$n_3$  &0&0\\ 
$\tau_4 $&  0& 0 &0&0\\  
$\tau_5$ & $2n_5$ &0  &0&0\\  
 \hline
\end{tabular}
\end{center}

The components of type $\tau_1$ and $\tau_2$ are the graphs of endomorphisms of the unit object $\unitob$ which are determined by their defect and excess.  The component of type $\tau_4$ may be ignored, since it gives the trivial endomorphism of $1$ as a morphism in $WW(\ilrel)$. Finally, the components of types 
$\tau_3$ and $\tau_5$ are determined by the shadow of $(f,g)$ with dimension $n_3$, and the space $X\times\zbar$ with dimension $2(n_3 + n_5)$.  It follows that the graph of our WW morphism, and thus the morphism itself, is determined by the shadow, the defect, and the excess.  

The constraints on the pair $(\cald,\cale)$ are just those from Proposition \ref{prop-defectexcessinequal}, while the composition law involving the stated cocycle follows again from  
Proposition
 \ref{prop-composedexcess} and Corollary \ref{cor-composeddefect}. 
\qed

\begin{rmk}
\label{rmk-free}
\emph{
Each of the endomorphism monoids above is free with two generators. The generators are $(0,1)$ and $(1,0)$ in the case of $\lrel$, $(1,1)$ and $(1,0)$ in the case of $\ilrel$, and $(0,1)$ and $(1,1)$ in the case of $\clrel$, 
}
\end{rmk}

\begin{rmk}
\label{rmk-presymplecticpoisson}
{\em
We recall that a (constant) presymplectic structure on $X$ is a possibly degenerate skew-symmetric bilinear form, while a (constant) Poisson structure on $X$ is a presymplectic structure on $X^*$.  The latter may be thought of as  a bracket operation on the linear functions on $X$ with values in the constants, leading to a (Heisenberg) Lie algebra structure on the affine functions.  The notions of isotropic and coisotropic subspace are natural in presymplectic and Poisson vector spaces respectively, giving rise to versions of $\ilrel$ and $\slrel$ involving these more general objects, and these larger categories are isomorphic via duality.  (The annihilator in $X^*$ of an isotropic subspace in presymplectic $X$ is coisotropic in Poisson $X^*$, and vice versa.)
 
The classification of (co)isotropic pairs in \cite{lo-we:coisotropic} should make it possible to extend the results here to these larger categories, but we have not yet attempted to carry out this extension.   The categories contain as monoidal subcategories both $\lrel$ and $\slrel$ (which generate all objects by monoidal product); the main difficulty is to understand how the morphisms interlace the two subcategories. }
\end{rmk}


\begin{thebibliography}{99}

\bibitem{be-tu:relazioni} Benenti, S., Tulczyjew, W., Relazioni lineari simplettiche, \emph{Memorie dell'Accademia delle Scienze di Torino} \textbf{5} (1981), 71-140.

\bibitem{br:endomorphisms} Brenner, S., Endomorphism algebras of vector spaces with distinguished sets of subspaces, \emph{Journal of Algebra} \textbf{6} (1967), 100-114.

\bibitem{et:introduction} Etingof, P., et al., \emph{Introduction to representation theory}, Lecture notes \url{www-math.mit.edu/~etingof/replect.pdf}. 

\bibitem{ge-po:problems} Gelfand, I. M., and Ponomarev, V. A., Problems of linear algebra and classification of quadruples of subspaces of a finite-dimensional vector space, \emph{Colloquia Mathematica Societatis Janos Bolyai} \textbf{5}, Tihany (1970), 163-237.

\bibitem{li-we:selective} 
Li-Bland, D., and Weinstein, A., Selective categories and linear canonical relations, \emph{SIGMA} \textbf{10} (2014), 100, 31 pages.

\bibitem{lo-we:coisotropic}
Lorand, J., and Weinstein, A., (Co)isotropic pairs in Poisson and presymplectic vector spaces, 
preprint arXiv 1503.00169.

\bibitem{ne:lectures}
Neretin, Yu. A., \emph{Lectures on Gaussian Integral Operators and Classical Groups}, European Mathematical Society, Z\"urich, 2011. 

\bibitem{to:linear} 
Towber, J., Linear relations, \emph{Journal of Algebra} \textbf{19} (1971), 1-20. 

\bibitem{we-wo:functoriality}
Wehrheim, K., and Woodward C.,
Functoriality for Lagrangian correspondences in Floer theory,
{\em Quantum topology} {\bf 1} (2010), 129--170.

\bibitem{we:note} 
Weinstein, A., A note on the Wehrheim-Woodward category, \emph{Journal of Geometric Mechanics} \textbf{3} (2011), 507-515. 

\end{thebibliography}
\end{document}